\documentclass[a4paper]{article}

\usepackage{amsmath, amsfonts, theorem, version, graphicx, lscape}
\usepackage{natbib}

\newtheorem{Prop}{Proposition}[section]

\newtheorem{rem}{Remark}
\def\Black{} 

\newcommand{\cqfd}{\hfill $\square$}

\newcommand{\R}{\mathbb R}

\newcommand{\n}{^{(n)}}

\newcommand{\Xb}{\mathbf{X}}
\newcommand{\Sb}{\mathbf{S}}

\newcommand{\vb}{\ensuremath{\mathbf{v}}}
\newcommand{\xb}{\ensuremath{\mathbf{x}}}

\newcommand{\Ub}{\ensuremath{\mathbf{U}}}

\newcommand{\tb}{\ensuremath{\mathbf{t}}}

\newcommand{\thetab}{{\pmb \theta}}

\newcommand{\varthetab}{{\pmb \vartheta}}

\newcommand{\Umb}{{\pmb \Upsilon}}

\newcommand{\zerob}{{\pmb 0}}
\newcommand{\Deltab}{{\pmb \Delta}}

\newcommand{\nub}{{\pmb \nu}}
\newcommand{\etab}{{\pmb \eta}}

\newcommand{\Gamb}{{\pmb \Gamma}}

\newcommand{\varepsb}{\pmb \varepsilon}

\newcommand{\pr}{^{\prime}}

\newcommand{\uut}[1] {\renewcommand{\arraystretch}{0.5}
\begin{array}[t]{c}{#1}\vspace{.2mm}\\
   \widetilde{}\vspace{-2.1mm}
\end{array}
\renewcommand{\arraystretch}{1} }

\newcommand{\utphi}{\!\!\! \uut{\phi}\!\!}
\newcommand{\utQ}{\!\! \uut{Q}\!\!}
\newcommand{\utDelta}{\!\! {\uut {\pmb \Delta}}\vspace{-1mm}\!\!}

\newcommand{\ny}{n\rightarrow\infty}

\begin{document}










%


\title{\bf Efficient ANOVA for directional data}

\author{{Christophe} {Ley} \thanks{
 D\'{e}partement de Math\'{e}matique and ECARES, Universit\'{e} Libre
 de Bruxelles 
Boulevard du Triomphe, CP 210 
B-1050 Bruxelles,
{Belgium}
}, {Yvik} {Swan}\thanks{Facult\'e des Sciences, de la Technologie et
  de la Communication, 
 Unit\'e de Recherche en Math\'ematiques, 
 Universit\'e de Luxembourg, 
 6, rue Richard Coudenhove-Kalergi, 
 L-1359 Luxembourg, {Grand Duch\'e de Luxembourg}
} and {Thomas} {Verdebout}\thanks{EQUIPPE, 
 Universit\' e Lille Nord de France, 
 Domaine Universitaire du Pont de Bois, BP 60149 
 F-59653 Villeneuve d'Ascq Cedex,  {France}
 }
}

\maketitle

\begin{abstract} 
  In this paper we tackle the ANOVA problem for directional data (with particular emphasis on geological data) by having recourse to the Le Cam methodology
  usually reserved for linear multivariate analysis.
  We construct locally and asymptotically most stringent parametric
  tests for ANOVA for directional data within the class of rotationally symmetric
  distributions. We turn these parametric tests into
  semi-parametric ones by (i) using a studentization argument (which
  leads to what we call pseudo-FvML tests) and by (ii) resorting to the
  invariance principle (which leads to efficient rank-based tests).  Within
  each construction the semi-parametric tests inherit optimality
  under a given distribution (the FvML distribution in the first case,
  any rotationally symmetric distribution in the second) from their
  parametric antecedents and also improve on the latter by being
  valid under the whole class of rotationally symmetric
  distributions. Asymptotic relative efficiencies are calculated and
  the finite-sample behavior of the proposed tests is investigated by
  means of a Monte Carlo simulation. We conclude by applying our findings on a real-data example involving geological data.

\


 \end{abstract}


\section{Introduction}\label{intro}

\setcounter{equation} {0}

Spherical or directional data naturally arise in a broad range of
earth sciences such as geology (see, e.g., Watson~1983 or Fisher and
Hall~1990), astrophysics, meteorology, oceanography or studies of
animal behavior (see, e.g., Merrifield~2006 and the references
provided therein) or even in neuroscience (see Leong and
Carlile~1998). Although primitive statistical analysis of directional
data can already be traced back to early 19th century works by the
likes of C.~F.~Gauss and D.~Bernoulli, the methodical and systematic
study of such non-linear data 
by means of tools tailored for their specificities  only begun in the
1950s under the impetus of Sir Ronald Fisher's pioneering work (see Fisher~1953). We refer the reader to the monographs
Fisher \emph{et al.}~(1987) and Mardia and Jupp~(2000) for a thorough 
introduction and comprehensive  overview   of this
discipline.


An important area of application of spherical statistics is in
geology (for instance for the study of palaeomagnetic data, see
McFadden and Jones~1981 or the more recent Acton~2011) wherein  the
data are usually   modeled as realizations of random vectors
$\Xb$ taking values on the surface of the unit hypersphere
$\mathcal{S}^{k-1}:=\{\vb\in\R^k:\vb'\vb=1\}$, 
the distribution of $\Xb$ depending only on its angular distance from
a fixed point $\thetab\in\mathcal{S}^{k-1}$ which is 
to be viewed as a ``north pole'' for the problem under study. 
%
%
A natural, flexible and realistic family of probability distributions
for such  data is the class of so-called \emph{rotationally symmetric distributions}
introduced by Saw~(1978) --  see  Section~\ref{rotsym} below for 
definitions and notations. Roughly speaking such distributions allow to
model all spherical  data that are spread out uniformly around a
central parameter $\thetab$ with the concentration of the data waning
as the angular distance from the north pole increases. 
Within this setup, an important    question  goes as follows : 
\emph{``do several measurements of remanent magnetization come from a
  same source of magnetism?''} More precisely, suppose
that there are $m$ different data sets 
spread around    $i$ sources of magnetism 
$\thetab_i \in\mathcal{S}^{k-1}$, $i = 
1,\ldots, m$.  The question then becomes that of testing for the
problem $\mathcal{H}_0: \thetab_1= \ldots= \thetab_{m}$ 
against $\mathcal{H}_1:\exists\, 1\leq i\neq j\leq m\,\,\mbox{such
  that}\,\,\thetab_i\neq\thetab_j$, that is, an ANOVA problem for directional data. 

This important  problem  has, obviously, already
been considered in the literature (see Mardia and Jupp~2000, chapter~10, for an overview). The difficulty of the task,
however, entails that most available methods are either of parametric nature or
suffer from computational difficulties/slowness such as
Wellner~(1979)'s permutation test or Beran and Fisher~(1998)'s
bootstrap test. To the best of our knowledge, the only computationally
simple and asymptotically distribution-free test for the general null
hypothesis $\mathcal{H}_0$ above is the test given in Watson~(1983).
%
The purpose of the present paper is to complement this literature by
constructing tests that are optimal under a given 
$m$-tuple of distributions--$(P_1,\ldots, P_m)$ say--but remain valid
(in the sense that they meet the nominal level constraint) under the general null hypothesis $\mathcal{H}_0$ involving a
large family of spherical distributions.  In particular,
 the tests we propose are asymptotically distribution-free within the semi-parametric
 class of rotationally symmetric distributions. 
Obviously the applicability of our ANOVA procedures is not reserved 
to geological data alone, but directly extends to any type of directional
data for which the assumption of  rotational symmetry with 
location parameter $\thetab$ seems to be reasonable. 

The backbone of our approach is the so-called Le Cam methodology (see
Le Cam~1986), as adapted to the spherical setup by Ley \emph{et
  al.}~(2013). Of utmost importance for our aims here is the
\emph{uniform local asymptotic normality} (ULAN) of a sequence of
rotationally symmetric distributions established therein and which we
adapt to our present purpose in Section~\ref{sec:param}. In the same
Section~\ref{sec:param} we also adapt results from Hallin \emph{et
  al.}~(2010) to determine the general form of a so-called
asymptotically \emph{most stringent} parametric test for the above
hypothesis scheme $\mathcal{H}_0$ 
against $\mathcal{H}_1$.
Due to its parametric nature the optimality of the $(P_1, \ldots,
P_m)$-parametric test is thwarted by its non-validity under any
$m$-tuple $(Q_1, \ldots, Q_m)$ distinct from $(P_1, \ldots, P_m)$.  In
order to palliate this problem we have recourse to two classical
tools which we adapt to the spherical setting: first a
\emph{studentization} argument, which leads to so-called
\emph{pseudo-Fisher-von Mises-Langevin} (pseudo-FvML) tests, and
second the \emph{invariance principle}, yielding optimal rank-based
tests. Both families of tests are of semi-parametric nature.

The idea behind the  {pseudo-FvML} test has the same flavor as the
\emph{pseudo-Gaussian} tests in the classical ``linear'' framework
(see, for instance, Muirhead and Waternaux~1980 or Hallin and
Paindaveine~2008 for more information on pseudo-Gaussian
procedures). More concretely, since the FvML distribution is
generally considered as the spherical analogue of the Gaussian distribution (see
Section~\ref{rotsym} for an explanation), our first approach consists
in  using the FvML as basis distribution 
and ``correcting'' the
(parametric) FvML most stringent test, optimal under a $m$-tuple
$(P_1, \ldots, P_m)$ of FvML distributions, in such a way that the
resulting test $\phi\n$ remains valid under the entire class of rotationally
symmetric distributions. We obtain the asymptotic distribution of the
asymptotically most stringent  {pseudo-FvML} test statistic $Q\n$
under the null and under contiguous alternatives. As it turns out, the test statistic
$Q\n$ and the test statistic $Q_{\rm Watson}\n$ provided in Watson
(1983) are asymptotically equivalent under the null (and therefore
under contiguous alternatives). As a direct consequence, we hereby
obtain, in passing, that Watson (1983)'s test  $\phi_{\rm
  Watson}\n$ also enjoys the property of being   asymptotically
most stringent in the FvML case.

The optimality property of $\phi_{\rm Watson}\n$ and $\phi\n$ is, by
construction,  restricted to situations in which the underlying
$m$-tuple of distributions is FvML. In the sequel we make use of the well-known invariance
principle to construct a more flexible
family of test statistics.   To this end we first obtain a group of monotone
transformations which generates the null hypothesis. Then we construct
tests based on the maximal invariant associated with this group. The
resulting tests (that are based on spherical signs and ranks) are,
similarly as $\phi_{\rm Watson}\n$ and $\phi\n$, asymptotically valid
under any $m$-tuple of rotationally symmetric densities. Our approach
here, however, further entails that for 
any given  $m$-tuple $(P_1, \ldots, P_m)$ of rotationally symmetric
distributions (not necessarily FvML ones) it suffices to choose the
appropriate  $m$-tuple
$\underline{K}=(K_1, \ldots, K_m)$ of score functions to guarantee
that  the resulting test is
asymptotically most stringent under $(P_1, \ldots, P_m)$.

The rest of the paper is organized as follows. In Section~\ref{rotsym}, we define
the class of rotationally symmetric distributions and collect
the main assumptions of the paper. In Section~\ref{sec:param}, we
summarize asymptotic results in the context of rotationally symmetric
distributions and show how to construct the announced optimal
parametric tests for the ANOVA problem. We then extend the latter to
pseudo-FvML tests in Section~\ref{sec:FvML} and to rank-based tests in
Section~\ref{sec:ranks}, and study their respective asymptotic
behavior in each section. Asymptotic relative efficiencies are
provided in Section~\ref{sec:ranks}. The theoretical results are
corroborated via a Monte Carlo simulation in
Section~\ref{sec:Simul}. A real data application is considered in
Section \ref{realdat}. Finally an appendix collects the proofs.

\section{Rotational symmetry}\label{rotsym}

Throughout, the $m(\geq2)$ samples of data points $\Xb_{i1}, \ldots, \Xb_{in_i}$, $i=1,\ldots,m$, are assumed to belong to the unit sphere $\mathcal{S}^{k-1}$ of $\R^k$, $k\geq2$, and to satisfy

\

\noindent {\sc Assumption A. (Rotational symmetry)} For all $i=1,\ldots,m$, $\Xb_{i1}, \ldots, \Xb_{in_i}$  are \mbox{i.i.d.} with common distribution ${\rm P}_{\thetab_i; f_i}$  characterized by a density (with respect to the usual surface area measure on spheres)
\begin{equation}\label{density} \xb\mapsto c_{k,f_i} \; f_i({\bf x}\pr\thetab_i),\quad\xb\in\mathcal{S}^{k-1},\end{equation} where $\thetab_i\in \mathcal{S}^{k-1}$ is a location parameter and $f_i: [-1,1]\rightarrow \R_0^+$ is absolutely continuous and (strictly) monotone increasing. Then,  if $\Xb$ has density \eqref{density},  the density of $\Xb\pr \thetab_i$ is of the form
$$t\mapsto\tilde{f}_i(t):= \frac{\omega_k \; c_{k,f_i}}{B(\frac{1}{2}, \frac{1}{2}(k-1))} f_i({t}) (1-t^2)^{(k-3)/2}, \quad -1\leq t \leq 1,$$
where $\omega_k= 2 \pi^{k/2}/ \Gamma (k/2)$ is the surface area of $\mathcal{S}^{k-1}$ and $B(\cdot,\cdot)$ is the beta function. The
corresponding cumulative distribution function (cdf) is denoted by $\tilde{F}_i(t)$, $i=1,\ldots,m$.

\

The functions $f_i$ are  called \emph{angular functions} (because the
distribution of each $\Xb_{ij}$ depends only on the angle between it
and the location $\thetab_i\in\mathcal{S}^{k-1}$). Throughout the rest
of this paper, we  denote by $\mathcal{F}^m$ the collection of $m$-tuples
of angular functions $\underline{f}:=(f_1,f_2,\ldots,f_m)$. Although not
necessary for the definition to make sense,  monotonicity of $f_i$
ensures that surface areas in the vicinity of the location parameter
$\thetab_i$ are allocated a higher probability mass than more remote
regions of the sphere. This property happens to be very appealing from
the modeling point of view. The assumption of rotational  symmetry
also entails  appealing stochastic properties.  Indeed, as shown
in Watson~(1983), for a random vector $\Xb$ distributed according to
some ${\rm P}_{\thetab_i; f_i}$ as in Assumption A,   not only is the multivariate sign
vector  ${\bf
  S}_{\thetab_i}(\Xb):=(\Xb-(\Xb'\thetab_i)\thetab_i)/||\Xb-(\Xb'\thetab_i)\thetab_i||$
uniformly distributed on ${\mathcal S}^{\thetab_i^\perp}:=\{ \vb
\in\R^k\,|\,\| \vb \|=1, \vb\pr\thetab_i=0 \}$ but also the angular
distance $\Xb\pr\thetab_i$ and the sign vector ${\bf S}_{\thetab_i}(\Xb)$
are stochastically independent.

The class of rotationally symmetric distributions contains a wide
variety of useful spherical distributions including  the   wrapped
normal distribution, the FvML, the linear, the logarithmic and the
logistic (a definition of the latter three is provided in Section~\ref{sec:ranks}
below).  The most popular and most
used rotationally symmetric distribution is the aforementioned  FvML distribution
(named, according to Watson~1983, after von Mises~1918, Fisher~1953,
and Langevin~1905), whose  density is of the form
\begin{equation*}\label{FvML}
f_{\rm FvML(\kappa)}(\xb;\thetab)=C_k(\kappa)\exp(\kappa\xb'\thetab),\quad \xb\in\mathcal{S}^{k-1}, 
\end{equation*}
where $\kappa>0$ is a concentration or dispersion parameter,
$\thetab\in\mathcal{S}^{k-1}$ a location parameter and  $C_k(\kappa)$
is the corresponding  normalizing constant.  
For ease of reference we shall, in what follows, rather use the
notation  $\phi_\kappa$ instead of  $f_{\rm
  FvML(\kappa)}$. This choice of notation is motivated both by the
wish for notational simplicity but also serves to  further underline
the analogy between  the FvML distribution as a spherical and the Gaussian distribution as a linear distribution. 
This  analogy is mainly due to the fact that
the FvML distribution is the only spherical distribution for which the
spherical empirical mean 
$\hat{\thetab}_{\rm
  Mean}:=\sum_{i=1}^n\Xb_i/||\sum_{i=1}^n\Xb_i||$
(based on observations $\Xb_1,\ldots,\Xb_n\in\mathcal{S}^{k-1}$) is  the Maximum Likelihood
Estimator (MLE) of its spherical location parameter, similarly as
  the Gaussian distribution  is the only (linear)
distribution in which the  empirical mean
$n^{-1}\sum_{i=1}^n\Xb_i$ (based on observations $\Xb_1,\ldots,\Xb_n\in\R^k$)  is  the MLE
for the (linear) location parameter.  
We refer the interested reader to  Breitenberger~(1963), Bingham and Mardia~(1975) or Duerinckx and
Ley~(2013) for details and references on this topic;  see also
Schaeben~(1992) for a discussion on spherical analogues of the  
Gaussian distribution. 


\section{ULAN and optimal parametric tests}\label{sec:param}

Throughout this paper a   test $\phi^*$ is called optimal if it is
\emph{most stringent} for testing $\mathcal{H}_0$ against
$\mathcal{H}_1$ within the class of tests $\mathcal{C}_\alpha$ of
level $\alpha$, that is if 
\begin{equation}
  \label{eq:1}
  \sup_{{\rm P}\in\mathcal{H}_1}r_{\phi^*}({\rm P})\leq\sup_{{\rm P}\in\mathcal{H}_1}r_\phi({\rm P})\quad\forall\phi\in\mathcal{C}_\alpha,
\end{equation}
where $r_{\phi_0}({\rm P})$ stands for the \emph{regret} of the test
$\phi_0$ under ${\rm P}\in\mathcal{H}_1$ defined as $r_{\phi_0}({\rm
  P}):=\left[\sup_{\phi\in\mathcal{C}_\alpha}{\rm E}_{\rm
    P}[\phi]\right]-{\rm E}_{\rm P}[\phi_0]$, the deficiency in power
of $\phi_0$ under ${\rm P}$ compared to the highest possible (for
tests belonging to $\mathcal{C}_\alpha$) power under ${\rm P}$.

As stated in the Introduction, the main ingredient for the
construction of optimal (in the sense of \eqref{eq:1}) parametric tests for the null hypothesis $\mathcal{H}_0: \thetab_1= \thetab_2= \ldots= \thetab_m$ consists in establishing the ULAN property of the parametric model $$\left(\left\{ {\rm       P}^{(n)}_{\thetab_1; f_1} \;  \vert \; \thetab_1 \in {\mathcal
      S}^{k-1} \right\},\ldots,\left\{ {\rm P}^{(n)}_{\thetab_m; f_m} \;
    \vert \; \thetab_m \in {\mathcal S}^{k-1} \right\}\right)$$  
for a
fixed $m$-tuple of (possibly different) angular functions $\underline{f}=(f_1,\ldots,f_m)$, where ${\rm
  P}^{(n)}_{\thetab_i; f_i}$ stands for the joint distribution of
$\Xb_{i1},\ldots,\Xb_{in_i}$, $i=1,\ldots,m,$ for a
fixed $m$-tuple of angular functions $(f_1, \ldots, f_m)$. Letting $\varthetab:=(\thetab_1\pr, \ldots, \thetab_m\pr)\pr$, we further denote by ${\rm
  P}^{(n)}_{\varthetab; \underline{f}}$ the joint law combining
${\rm P}^{(n)}_{\thetab_1; f_1}$, \ldots, ${\rm P}^{(n)}_{\thetab_m;
  f_m}$. In order to be able to state our results, we need to impose a
certain amount of control on the respective sample sizes $n_i$, $i=1, \ldots,
m$. This we   achieve via the following 

\

\noindent {\sc Assumption B}. Letting $n=\sum_{i=1}^m n_i$, for all $i=1, \ldots, m$ the ratio
$r_i\n:=n_i/n$ converges to a finite constant $r_i$ as $\ny$. 

\

\noindent In particular Assumption~B entails that the specific
sizes $n_i$ are, up to a point, irrelevant; hence  in what precedes and in what
follows, we simply use the superscript $^{(n)}$  for the different
quantities at play and do not specify whether they are associated with
a given $n_i$.  In the sequel we let ${\rm diag}({\bf A}_1, \ldots, {\bf A}_m)$ stand
for the $m \times m$ block-diagonal matrix with blocks ${\bf A}_1,
\ldots, {\bf A}_m$,  and use the notation 
$  \nub\n:={\rm
  diag}((r_1\n)^{-1/2} {\bf I}_k, \ldots, (r_m\n)^{-1/2} {\bf I}_k).
$

  Informally, a sequence of rotationally symmetric models $\left\{ {\rm
  P}^{(n)}_{\varthetab; \underline{f}} \;  \vert \;
  \varthetab \in ({\mathcal S}^{k-1})^m \right\}$ is ULAN if,
uniformly in  $\varthetab\n=(\thetab_1^{(n) \prime}, \ldots, \thetab_m^{(n)
  \prime})\pr \in(\mathcal{S}^{k-1})^m$ such that 
$\varthetab\n-\varthetab=O(n^{-1/2})$,   the log-likelihood 
\begin{equation*}
 \log  \left({\rm P}^{(n)}_{\varthetab\n+n^{-1/2} \nub\n {\tb}\n;
  \underline{f}}/{\rm P}^{(n)}_{\varthetab\n;
  \underline{f}}\right)
\end{equation*}
allows a specific form of (probabilistic) Taylor expansion (see
equation \eqref{eq:2} below) as a function of 
$\tb\n:=(\tb_1^{(n)   \prime},\ldots, \tb_m^{(n) \prime})\pr \in
\R^{mk}$.  Of course the local perturbations $\tb\n$ must be chosen so
that  $\varthetab\n+ n^{-1/2} \nub\n \tb^{(n)}$ remains
on $(\mathcal{S}^{k-1})^m$ and thus, in particular, the  $\tb_i\n$
need to satisfy
\begin{align} \label{localt} 
0 &= (\thetab_i\n+ n_i^{-1/2}\tb_i\n)\pr ({\thetab_i\n}+ n_i^{-1/2}\tb_i\n) -1 \nonumber \\ 
&=  2n_i^{-1/2}({\thetab_i\n})\pr \tb_i\n + n_i^{-1} (\tb_i\n)\pr
\tb_i\n \end{align}
for all $i=1, \ldots,m$.
Consequently, $\tb_i\n$ must be such that 
$2 n_i^{-1/2}({\thetab_i\n})\pr   \tb_i\n+
o(n_i^{-1/2})=0$ : for ${\thetab_i\n}+
n_i^{-1/2}\tb_i\n$ to remain in ${\mathcal S}^{k-1}$, the perturbation
$\tb_i\n$ must
belong, up to a $o(n_i^{-1/2})$ quantity, to the tangent space to
$\mathcal{S}^{k-1}$ at $\thetab_i\n$.

The domain of the parameter being the non-linear manifold $\left({\mathcal
  S}^{k-1}\right)^m$ it is all but easy to establish the ULAN property of a
sequence of rotationally symmetric models. A natural way to 
handle this difficulty consists, as in  Ley \emph{et al.}~(2013),  in resorting to a  re-parameterization
of the problem in terms of spherical coordinates $\etab$, say, for which it is possible to
prove ULAN, subject to the following technical condition on the
angular functions. 

\

\noindent {\sc Assumption C}. The Fisher information associated with the
spherical location parameter is finite; this finiteness is ensured if,
for $i=1, \ldots, m$ and letting $\varphi_{f_i}:=\dot{f}_i/f_i$ ($\dot{f_i}$ is
the a.e.-derivative of $f_i$), $\mathcal{J}_k({f_i}):= \int_{-1}^{1}
\varphi_{f_i}^2(t) (1-t^2) \tilde{f}_i(t) dt < + \infty$.

\

\noindent After obtaining the ULAN property for the 
$\eta$-parameterization, one can use a lemma from Hallin \emph{et
  al.}~(2010) to transpose the ULAN property in the spherical
$\etab$-coordinates back in terms of the  original
$\thetab$-coordinates.  Finally the inner-sample
independence and the mutual independence between the $m$ samples
entail that we can deduce the required   ULAN property which is
relevant for our purposes (this  we state without proof because it follows directly  from
Proposition~2.2 of  Ley \emph{et al.} 2013).

\begin{Prop}\label{ULAN} Let Assumptions A, B and C hold. Then the model $\left\{ {\rm P}\n_{\varthetab; \underline{f}} \;  \vert \; \varthetab \in ({\mathcal S}^{k-1})^m \right\}$ is ULAN with central sequence  $\Deltab\n_{\varthetab ; \underline{f}}:=\left((\Deltab^{(n)}_{\thetab_1; f_1 })',\ldots, (\Deltab^{(n)}_{\thetab_m;f_m })' \right)\pr$, where
\begin{equation*}\label{centrseqm}\Deltab\n_{\thetab_i; f_i}:= n_i^{-1/2} \sum_{j=1}^{n_i} \varphi_{f_i}(\Xb_{ij} \pr \thetab_i) (1-(\Xb_{ij} \pr \thetab_i)^2)^{1/2} \Sb_{\thetab_i}(\Xb_{ij}),\quad i=1,\ldots,m, \end{equation*}
and Fisher information matrix $\Gamb_{\varthetab;\underline{f}}:={\rm diag}(\Gamb_{\thetab_1; f_1},\ldots,\Gamb_{\thetab_m; f_m})$
where
\begin{equation*} \label{Fishm} \Gamb_{\thetab_i; f_i}:= \frac{\mathcal{J}_k({f_i})}{k-1} ({\bf I}_k- \thetab_i \thetab_i\pr),\quad i=1,\ldots,m. \end{equation*}
More precisely, for any $\varthetab\n \in(\mathcal{S}^{k-1})^m$ such that $\varthetab\n-\varthetab=O(n^{-1/2})$ and any bounded sequences $\tb\n=(\tb^{(n) \prime}_1,\ldots,\tb^{(n) \prime}_m)\pr$ as in~(\ref{localt}), we have
\begin{equation}
  \label{eq:2}
  \log \left( \frac{{\rm P}\n_{\varthetab\n + n^{-1/2} \nub\n \tb\n ; \underline{f}}}{{\rm P}\n_{\varthetab\n; \underline{f}}}\right)=   (\tb\n)'\Deltab\n_{\varthetab\n; \underline{f}} - \frac{1}{2}(\tb\n)' \Gamb_{\varthetab ; \underline{f}} \tb\n+o_{\rm P}(1),
\end{equation}
where $\Deltab\n_{\varthetab\n; \underline{f}}\stackrel{\mathcal{L}}{\rightarrow}\mathcal{N}_{mk}(\zerob,\Gamb_{\varthetab; \underline{f}})$, both under ${\rm P}\n_{\varthetab; \underline{f}}$, as $\ny$.
\end{Prop}

Proposition \ref{ULAN} provides us with all the necessary tools for
building optimal $\underline{f}$-parametric procedures (i.e. under any
$m$-tuple of densities with respective specified angular functions
$f_1, \ldots, f_m$) for testing $\mathcal{H}_0:\thetab_1= \ldots=
\thetab_m$ against $\mathcal{H}_1:\exists\, 1\leq i\neq j\leq
m\,\,\mbox{such that}\,\,\thetab_i\neq\thetab_j$.  Intuitively, this
follows from the fact that the second-order expansion of the
log-likelihood ratio for the model $\left\{ {\rm P}\n_{\varthetab;
    \underline{f}} \; \vert \; \varthetab \in ({\mathcal S}^{k-1})^m
\right\}$ strongly resembles the log-likelihood ratio for the
classical Gaussian shift experiment, for which optimal procedures are
well-known and are based on the corresponding first-order term. 
Now clearly the null hypothesis $\mathcal{H}_0$ is the intersection
between $(\mathcal{S}^{k-1})^m$ and the linear subspace (of
$\R^{mk}$) $${\mathcal C} :=\{ \vb=(\vb_1\pr, \ldots, \vb_m \pr
)\pr\,|\, \vb_1, \ldots,\vb_m \in \R^k\, {\rm and} \;\vb_1=\ldots=
\vb_m\} =: {\mathcal M}({\bf 1}_m \otimes {\bf I}_k)$$ where we put
${\bf 1}_m:=(1, \ldots, 1)\pr \in \R^m$, ${\mathcal M}({\bf A})$ for the
linear subspace spanned by the columns of the matrix ${\bf A}$ and
${\bf A} \otimes {\bf B}$ for the Kronecker product between
${\bf A}$ and ${\bf B}$. Such a restriction, namely an intersection
between a linear subspace and a non-linear manifold, has already been
considered in Hallin \emph{et al.} (2010) in the context of Principal
Component Analysis (in that paper, the authors obtained very general
results related to hypothesis testing in ULAN families with curved
experiments). In particular from their results we can deduce that,
in order to obtain a locally and asymptotically most stringent test
in the present context, one has to consider the locally and
asymptotically most stringent test for the (linear) null hypothesis
defined \Black by the intersection between $\mathcal{C}$ and the
tangent to $(\mathcal{S}^{k-1})^m$. Let $\thetab$ denote the common
value of $\thetab_1, \ldots, \thetab_m$ under the null. In the
vicinity of ${\bf 1}_m \otimes \thetab$, the intersection between
$\mathcal{C}$ and the tangent to $(\mathcal{S}^{k-1})^m$ is given by
\begin{eqnarray} \label{syst} & \hspace{-40mm} \Big\{ (\thetab\pr + n^{-1/2} (r_1\n)^{-1/2} \tb_1^{(n)\prime}, \ldots,  \thetab\pr + n^{-1/2} (r_m\n)^{-1/2} \tb_m^{(n)\prime})\pr, 
\\ \hspace{20mm} & \hspace{20mm} {\thetab}\pr\tb_1\n=\ldots
={\thetab}\pr\tb_m\n=0, (r_1\n)^{-1/2} \tb_1\n= \ldots= (r_m\n)^{-1/2}
\tb_m\n \Big\}. \nonumber  \end{eqnarray} 
Solving the system (\ref{syst}) yields \begin{equation} \label{linres}
  \nub\n \tb\n=\left((r_1\n)^{-1/2} {\bf t}_1^{(n)\prime}, \ldots,
    (r_m\n)^{-1/2} \tb_m^{(n)\prime}\right)\pr \in {\mathcal M}({\bf
    1}_{m} \otimes ({\bf I}_k- \thetab\thetab\pr)). \end{equation} 
Loosely speaking we have ``transcripted'' the initial null
hypothesis $\mathcal{H}_0$ into a linear restriction of the form
(\ref{linres}) in terms of local perturbations $\tb\n$, for which Le
Cam's asymptotic theory  then  provides a locally and asymptotically
optimal parametric test under fixed $\underline{f}$. Using
Proposition \ref{ULAN} and letting $\Umb_{\varthetab}:={\bf 1}_{m}
\otimes ({\bf I}_k- \thetab\thetab\pr)$ and $\Umb_{\varthetab;
  {\nub}}\n:=(\nub\n)^{-1}\Umb_{\varthetab}$, an asymptotically  {most
  stringent} test $\phi_{\underline{f}}$ is then obtained by rejecting
$\mathcal{H}_0$ as soon as (${\bf A}^-$ stands for the Moore-Penrose
pseudo-inverse of ${\bf A}$) 
\begin{equation}\label{paramtest}Q_{\underline{f}}\n:=
  \Deltab_{\varthetab; \underline{f}}\pr \left(\Gamb_{\varthetab;
      \underline{f}}^{-} - \Umb_{\varthetab; {\nub}}\n
    \left((\Umb_{\varthetab; {\nub}}\n)\pr\Gamb_{ \varthetab;
        \underline{f}} \Umb_{\varthetab; {\nub}}\n \right)^{-}
    (\Umb_{\varthetab; {\nub}}\n)\pr\right) \Deltab_{\varthetab;
    \underline{f}} \end{equation} 
exceeds the $\alpha$-upper quantile of a chi-square distribution with
$(m-1)(k-1)$ degrees of freedom.  Hence the optimal parametric tests
are now known. 

There nevertheless remains much work to do. Indeed 
not only does the optimality  of our test $\phi_{\underline{f}}$   only hold  under the
$m$-tuple of angular densities $\underline{f}=(f_1, \ldots, f_m)$, but
also this parametric test suffers from the (severe) drawback of being
only valid under that pre-specified $m$-tuple. Since it is highly
unrealistic in practice to assume that 
the underlying densities are known, these tests are useless for
practitioners. Moreover, we so far have assumed known the common value of the spherical location under the null, which is unrealistic, too. The next two sections contain two distinct solutions
allowing to set these problems right.

\section{Pseudo-FvML tests}\label{sec:FvML}

For a given $m$-tuple of FvML densities $(\phi_{\kappa_1}, \ldots, \phi_{\kappa_m})$
with respective concentration parameters $\kappa_1, \ldots, \kappa_m>0$ (where
we do not assume $\kappa_1=\ldots= \kappa_m$), the score functions 
$\varphi_{\phi_{\kappa_i}}$ reduce to the constants $\kappa_i$,
$i=1, \ldots, m$, and hence the central sequences for each sample take
the simplified form
\begin{eqnarray*}
\label{centrseqfvml}
\Deltab\n_{\thetab_i; \phi_{\kappa_i} }&:=&\kappa_i n_i^{-1/2}
\sum_{j=1}^{n_i}  (1-(\Xb_{ij} \pr \thetab_i)^2)^{1/2}
\Sb_{\thetab_i}(\Xb_{ij}) \nonumber \\ &=& \kappa_i n_i^{-1/2}
\sum_{j=1}^{n_i}  ({\Xb}_{ij} -(\Xb_{ij}\pr \thetab_i) \thetab_i) \\ &
= & \kappa_i ({\bf I}_k- \thetab_i\thetab_i\pr) n_i^{-1/2}
\sum_{j=1}^{n_i} \Xb_{ij} \\ &=:& \kappa_i ({\bf I}_k-
\thetab_i\thetab_i\pr) n_i^{1/2} \bar{\Xb}_i\\&=& \kappa_i ({\bf I}_k-
\thetab_i\thetab_i\pr) n_i^{1/2} (\bar{\Xb}_i- \thetab_i),\quad i=1, \ldots, m.
\end{eqnarray*}
Optimal FvML-based procedures (in the sense of \eqref{eq:1}) for ${\mathcal H}_0$ are then built upon $\Deltab\n_{\varthetab;
 \underline{\phi}}=(\Deltab^{(n) \prime}_{\thetab_1; \phi_{\kappa_1}
}, \ldots, \Deltab^{(n) \prime}_{\thetab_m; \phi_{\kappa_m} })\pr$,
where $\underline{\phi}:=(\phi_{\kappa_1}, \ldots, \phi_{\kappa_m})$.

Before proceeding we here again draw the reader's
attention to the fact that  a parametric test built upon $\Deltab\n_{\varthetab;
 \underline{\phi}}$ will  only be  valid under
the $m$-tuple $\underline{\phi}$ and becomes non-valid
even if only the concentration parameters change. In this section,
this non-validity problem will be overcome in the following way. We
will first study the asymptotic behavior of
$\Deltab\n_{\varthetab;
 \underline{\phi}}$
under any given $m$-tuple $\underline{g}=(g_1, \ldots, g_m)\in\mathcal{F}^m$ and
consider the newly obtained quadratic form in
$\Deltab\n_{\varthetab;
 \underline{\phi}}$. Clearly, this quadratic form
will now depend on the asymptotic variance of
$\Deltab\n_{\varthetab;
 \underline{\phi}}$
under $\underline{g}$, hence again, for each $\underline{g}$, we are
confronted to an only-for-$\underline{g}$-valid test statistic. The
next  step then consists in applying a studentization
argument, meaning that we replace the asymptotic variance quantity by
an appropriate estimator. We then study the asymptotic behavior of the
new quadratic form under any 
$m$-tuple of rotationally symmetric distributions. As we will show, the final outcome of
this procedure will be tests which happen to be optimal under
\emph{any} $m$-tuple of FvML distributions (that is, for any values
$\kappa_1, \ldots, \kappa_m>0$) and valid under the entire class of
rotationally symmetric distributions; these tests are our so-called
pseudo-FvML tests.

For the sake of readability, we adopt in the sequel the notations ${\rm
  E}_f[\cdot]$ for expectation under the angular function $f$ and $\varthetab_0:=:{\bf 1}_m \otimes \thetab$ (where we recall that $\thetab$ represents the common value of $\thetab_1,\ldots,\thetab_m$ under the null). The
following result characterizes, for a given $m$-tuple of angular functions
$\underline{g}\in\mathcal{F}^m$, the asymptotic properties of the
FvML-based central sequence $\Deltab\n_{\varthetab_0;
 \underline{\phi}}$, both under ${\rm P}\n_{\varthetab_0;
 \underline{g}}$ and ${\rm P}\n_{\varthetab_0 +n^{-1/2} \nub\n \tb\n;
 \underline{g}}$ with $\tb^{(n)}$ as in~(\ref{localt}) for each sample.

\begin{Prop}\label{lineasympbis}
Let Assumptions A, B and C hold. Then, letting $B_{k, g_i}:=1-{\rm E}_{g_i}\left[(\Xb_{ij}\pr\thetab)^2\right]$ for $i=1, \ldots,m$, we have that  $\Deltab\n_{\varthetab_0; \underline{\phi}}$ is  
\begin{itemize}
\item [(i)] asymptotically normal under ${\rm P}\n_{\varthetab_0; \underline{g}}$ with mean zero and covariance matrix 
$$\Gamb^*_{\varthetab_0; \underline{g}}:={\rm diag}\left(\Gamb^*_{\thetab;
    g_1}, \ldots  \Gamb^*_{\thetab;  g_m} \right),$$ 
where $ \Gamb^*_{\thetab;  g_i}:=\frac{\kappa_{i}^2 B_{k, g_i}}{k-1} ({\bf I}_k -\thetab \thetab\pr), \quad i=1, \ldots, m;$
\item [(ii)] asymptotically normal under ${\rm
    P}\n_{\varthetab_0+n^{-1/2} \nub\n \tb\n; \underline{g}}$ ($\tb^{(n)}$ as
  in~(\ref{localt})) with mean $ \Gamb_{\varthetab_0; \underline{\phi},
    \underline{g}} \tb$ ($\tb:=(\tb_1', \ldots, \tb_m')'$ with $\tb_i:=
  \lim_{\ny} \tb_i^{(n)}$, $i=1, \ldots, m$) and
  covariance matrix $\Gamb^*_{\varthetab_0; \underline{g}}$, where, putting
  $C_{k,g_i}:={\rm E}_{g_i}[(1-(\Xb_{i j}\pr \thetab)^2)
  \varphi_{g_i}(\Xb_{ij}\pr \thetab)]$ for $i=1,\ldots, m$, 
$$\Gamb_{\varthetab_0; \underline{\phi}, \underline{g}} := {\rm diag}\left(\Gamb_{\thetab; \phi_{\kappa_1},  g_1},  \ldots, \Gamb_{\thetab;  \phi_{\kappa_m}, g_m} \right)$$ with
 $\Gamb_{\thetab; \phi_{\kappa_i},  g_i}:=\frac{\kappa_{i} C_{k,g_i}}{k-1} ({\bf I}_k -\thetab \thetab\pr),\quad i=1, \ldots, m.$
\end{itemize}
\end{Prop}

See the Appendix for the proof. As the null hypothesis only specifies
that the spherical locations coincide, we need to estimate the unknown
common value $\thetab$. Therefore, we assume the
existence of an estimator $\hat{\thetab}$ of $\thetab$ such that the
following assumption holds.  

\

\noindent {\sc Assumption D}.
 The estimator $\hat{\varthetab}={\bf 1}_m \otimes \hat{\thetab}$, with $\hat{\thetab} \in {\mathcal S}^{k-1}$, is \textit{$n^{1/2}\big(\nub\n\big)^{-1}$-consistent:}  for all $\varthetab_0={\bf 1}_m \otimes \thetab \in {\mathcal H}_{0}$,  $n^{1/2} \big(\nub\n\big)^{-1}(\hat{\varthetab} -\varthetab_0)=O_{\rm P}(1)$\Black, as $n\rightarrow \infty$ under ${\rm P}\n_{\varthetab_0 ; \underline{g}}$
 for any $\underline{g} \in{\mathcal F}^m$.  

\
  
  \noindent Typical examples of estimators satisfying Assumption D
  belong to the class of $M$-estimators (see Chang~2004) or
  $R$-estimators (see Ley \emph{et al.}~2013).  Put simply, instead of
  $\Deltab\n_{\varthetab_0; \underline{\phi}}$ we have to work with
  $\Deltab\n_{\hat\varthetab; \underline{\phi}}$ for some estimator
  $\hat{\varthetab}$ satisfying Assumption D. The next crucial result
 quantifies in how far this replacement affects the asymptotic
  properties established in Proposition~\ref{lineasympbis} (a proof is provided in the Appendix).

\begin{Prop}\label{cross}
Let Assumptions A, B and C hold and let $\hat{\varthetab}={\bf 1}_m \otimes \hat\thetab$ be an estimator of  $\varthetab_0$ such that Assumption D holds. Then
\begin{itemize}
\item[(i)] letting $\Umb\n:= \left(\sqrt{r_1\n}\; {\bf I}_k \;   \vdots  \ldots \vdots \;  \sqrt{r_m\n} \; {\bf I}_k\right)\pr$, $\Deltab\n_{\varthetab_0; \underline{\phi}}$ satisfies, under ${\rm P}\n_{\varthetab_0; \underline{g} }$ and as $\ny$,
$$\Deltab\n_{\hat\varthetab; \underline{\phi}}- \Deltab\n_{\varthetab_0; \underline{\phi}} =-\Gamb_{\varthetab_0; \underline{g}}^{\underline{\phi}} \Umb\n \sqrt{n} \left(\hat\thetab -\thetab\right)+ o_{\rm P}(1),$$
where $$\Gamb_{\varthetab_0;  \underline{g}}^{\underline{\phi}} :={\rm
  diag} \left(\Gamb_{\thetab; g_1}^{\phi_{\kappa_1}}, \ldots, \Gamb_{\thetab;
    g_m}^{\phi_{\kappa_m}} \right)$$ 
with $\Gamb_{\thetab; g_i}^{\phi_{\kappa_i}}:=\kappa_{i}{\rm
  E}_{g_i}\left[\Xb_{ij}\pr\thetab\right] ({\bf I}_k -\thetab
\thetab\pr), \quad i=1, \ldots, m;$ 

\item[ii)] for all $\varthetab \in ({\mathcal S}^{k-1})^m$,
  $\Gamb_{\varthetab;\underline{\phi},\underline{\phi}}=\Gamb_{\varthetab;\underline{\phi}}^{\underline{\phi}}=
  \Gamb_{\varthetab; \underline{\phi}}^*.$   

\end{itemize}
\end{Prop}

Following the inspiration of  Hallin and Paindaveine~(2008) (where a very general
theory for pseudo-Gaussian procedures is described) we are in a
position to  use 
Proposition  \ref{cross} to construct our pseudo-FvML tests. To this
end define, for $i=1, \ldots,
m$,   the quantities ${E}_{k,
  g_i}:={\rm E}_{g_i}\left[\Xb_{i j}\pr\thetab\right] $,  and set, for notational simplicity, ${D}_{k, g_i}:={E}_{k, g_i}/{B}_{k,
  g_i}$  and
${H}_{\underline{\phi},\underline{g}}:=\sum_{i=1}^m r_i\n {D}_{k,
  g_i}^2 {B}_{k, g_i}$. Then, letting 
\begin{align*} & \Gamb_{\varthetab_0; \underline{\phi},
    {\underline g}}^{\perp} := (\Gamb_{\varthetab_0; {\underline g}}^*)^-\nonumber \\
&  \quad  -    \!(\Gamb_{\varthetab_0; {\underline g}}^*)^- \Gamb_{\varthetab_0; {\underline g}}^{{\underline\phi}} \Umb_{\varthetab_0; {\nub}}\n  [(\Umb_{\varthetab_0; {\nub}}\n)\pr \Gamb_{\varthetab_0; {\underline g}}^{{\underline\phi}}  (\Gamb_{\varthetab_0; {\underline g}}^*)^-  \Gamb_{\varthetab_0; {\underline g}}^{{\underline \phi}} \Umb_{\varthetab_0; {\nub}}\n]^- (\Umb_{\varthetab_0; {\nub}}\n)\pr\Gamb_{\varthetab_0; {\underline g}}^{{\underline \phi}}  (\Gamb_{\varthetab_0; {\underline g}}^*)^- , \end{align*}
the $\underline{g}$-valid test statistic for $\mathcal{H}_0: \thetab_1= \ldots =\thetab_m$ we propose is the quadratic form
$$Q\n(\underline{g}):=(\Deltab\n_{\hat{\varthetab};
  \underline{\phi}})\pr {\Gamb}_{\hat{\varthetab}; \underline{\phi},
  \underline{g}}^{\perp}\Deltab\n_{\hat{\varthetab};
  \underline{\phi}}$$
It is
easy to verify that  $Q\n(\underline{g})$ does not
depend explicitly on the underlying concentrations $\kappa_1, \ldots,
\kappa_m$ but still depends on the  quantities $B_{k, g_i}$ and
$E_{k,g_i}$, $i=1,\ldots, m$. This obviously  hampers the validity of the statistic
outside of $\underline{g}$. The last step thus consists in estimating
these quantities. Consistent (via the Law of Large Numbers) estimators
for each of them are provided by $\hat{B}_{k, g_i}:=1-n_i^{-1}
\sum_{j=1}^{n_i} (\Xb_{i j}\pr \hat{\thetab})^2$ and $\hat{E}_{k,
  g_i}:= n_i^{-1} \sum_{j=1}^{n_i} (\Xb_{i j}\pr \hat{\thetab})$,
$i=1, \ldots, m$. For the sake of readability, we naturally also use
the notations $\hat{D}_{k, g_i}:=\hat{E}_{k, g_i}/\hat{B}_{k, g_i}$,
$i=1, \ldots, m$, and
$\hat{H}_{\underline{\phi},\underline{g}}:=\sum_{i=1}^m r_i\n
\hat{D}_{k, g_i}^2 \hat{B}_{k, g_i}$. Putting
$\bar{\Xb}_i:=n_i^{-1}\sum_{j=1}^{n_i}\Xb_{ij}$ for all
$i=1,\ldots,m$, straightforward calculations then show that our
pseudo-FvML test statistic for the $m$-sample spherical location
problem is  
\begin{align*}
 Q\n =& (k-1) \sum_{i=1}^{m}  \frac{n_i\hat{D}_{k, g_i}}{\hat{E}_{k, g_i}} \bar{\Xb}_i\pr ({\bf I}_k -{\hat\thetab} {\hat \thetab}\pr) \bar{\Xb}_i -(k-1)\sum_{i, j}^{m} \frac{n_i n_{j}}{n} \frac{\hat{D}_{k,g_i}\hat{D}_{k,g_{j}}}{\hat{H}_{\underline{\phi},\underline{g}}}  \bar{\Xb}_i\pr ({\bf I}_k -{\hat \thetab} {\hat \thetab}\pr) \bar{\Xb}_{j}, 
\end{align*}
which no more depends on $\underline{g}$. 

The following proposition, whose proof is given in the Appendix,
finally yields the asymptotic properties of this quadratic form under
the entire class of rotationally symmetric distributions, showing that
the test is well valid under that broad set of distributions. 

\begin{Prop}\label{chideuxtwo}
Let Assumptions A, B and C hold and let $\hat{\varthetab}$ be an estimator of $\varthetab_0$ such that Assumption D holds. Then
\begin{itemize}
\item[(i)] $Q\n$ is asymptotically chi-square with $(m-1)(k-1)$ degrees of freedom under \linebreak  $\bigcup_{\varthetab_0\in\mathcal{H}_0}\bigcup_{\underline{g} \in {\mathcal F}^m} {\rm P}\n_{{\varthetab}_0; \underline{g}}$;
\item[(ii)] $Q\n$ is asymptotically non-central chi-square with $k-1$ degrees of freedom and non-centrality parameter 
$$ l_{{\varthetab}_0,{\bf t}; \underline{\phi}, {\underline{g}}}:=\tb\pr \Gamb_{\varthetab_0; \underline{\phi}, \underline{g}} {\Gamb}_{{\varthetab}_0; \underline{\phi},\underline{g}}^{\perp} \Gamb_{{\varthetab}_0; \underline{\phi}, \underline{g}} \tb$$ under ${\rm P}\n_{{\varthetab}_0+n^{-1/2} \nub\n \tb\n; \underline{g}}$, where $\tb^{(n)}$ is as in~(\ref{localt}) and  $\tb:= \lim_{\ny} \tb^{(n)}$;
 \item[(iii)] the test $\phi\n$ which rejects the null hypothesis as soon as $Q\n$ exceeds the $\alpha$-upper quantile of the chi-square distribution with $(m-1)(k-1)$ degrees of freedom has 
asymptotic level $\alpha$ under $\bigcup_{\varthetab_0\in {\mathcal H}_{0}}\bigcup_{{\underline g} \in\mathcal{F}^{m}} \{ {\rm P}\n _{\varthetab_0; {\underline g}} \}$;
\item[(iv)]$\phi\n$ is locally and asymptotically most stringent, at asymptotic level~$\alpha$,  \linebreak for  $\bigcup_{\varthetab_0 \in {\mathcal H}_{0}} \bigcup_{{\underline g} \in\mathcal{F}^{m}}
\{ {\rm P}\n _{\varthetab_0; {\underline g}} \}$ against alternatives of the form
 $\bigcup_{\varthetab \notin{\mathcal H}_{0}} \{ {\rm
  P}\n _{\varthetab;{\underline \phi}}\}$.
  \vspace{1mm}
\end{itemize}
\end{Prop}

\begin{rem}
  It is easy to verify that $Q\n$ is asymptotically equivalent (the difference is
a $o_{\rm P}(1)$ quantity) to the test statistic for the same problem proposed in
Watson (1983) under the null (and therefore also under contiguous
alternatives). Thus, although the construction we propose is
different,  our pseudo-FvML tests coincide with Watson's proposal. In
passing, we have therefore also proved the asymptotic most stringency
of the latter.
\end{rem}

\section{Rank-based tests}\label{sec:ranks}

The pseudo-FvML test constructed in the previous section is valid
under any $m$-tuple of (non-necessarily equal) rotationally symmetric
distributions and retains the optimality properties of the FvML
most stringent parametric test in the FvML case. Although the FvML assumption
is often reasonable in practice, our aim in the present
section is to depart from this assumption and provide tests which are
optimal under any distribution. 

We start from
any given $m$-tuple $\underline{f}\in\mathcal{F}^m$ and our objective is to
turn the $\underline{f}$-parametric tests into tests which are still
valid under any $m$-tuple of (non-necessarily equal) rotationally symmetric
distributions and which remain optimal under $\underline{f}$. To
obtain such a test, we have recourse here to the second of the
aforementioned tools to turn our parametric tests into semi-parametric
ones: the invariance principle. This principle advocates that, if the
sub-model identified by the null hypothesis is invariant under the
action of a group of transformations $\mathcal{G}_T$, one should
exclusively use procedures whose outcome does not change along the
orbits of that group $\mathcal{G}_T$. This is the case if and only if
these procedures are measurable with respect to the maximal invariant
associated with $\mathcal{G}_T$. The invariance principle is
accompanied by an appealing corollary for our purposes here:
provided that the group $\mathcal{G}_T$ is a generating group for
$\mathcal{H}_0$, the invariant procedures are distribution-free under
the null.

Invariance with respect to ``common rotations" is crucial in this
context. More precisely, letting ${\bf O}\in {\mathcal SO}_k:=\{{\bf
  A} \in \R^{k \times k}, {\bf A}\pr{\bf A}={\bf I}_k, \;{\rm
  det}({\bf A})=1\}$, the null hypothesis is unquestionably invariant
with respect to a transformation of the form
$$g_{\bf O}: \Xb_{11}, \ldots, \Xb_{1n_1}, \ldots, \Xb_{m1}, \ldots, \Xb_{m
  n_{m}} \mapsto  {\bf O}\Xb_{11}, \ldots, {\bf O}\Xb_{1n_1},  \ldots, {\bf
  O}\Xb_{m1}, \ldots, {\bf O}\Xb_{m n_{m}}.$$
However, this group is not a generating for $\mathcal{H}_0$ as it does
not take into account the underlying angular functions
$\underline{f}$, which are an infinite-dimensional nuisance under
$\mathcal{H}_0$. This group is actually rather generating for
$\bigcup_{\varthetab_0\in\mathcal{H}_0}{\rm
  P}\n_{\varthetab_0;\underline{f}}$ with fixed
$\underline{f}$. Now, denote as in the previous section  the common
value of $\thetab_1, \ldots, \thetab_m$ under the null as $\thetab$. Then $\Xb_{ij} = (\Xb_{ij}\pr \thetab)  \thetab +
\sqrt{1- (\Xb_{ij}\pr\thetab)^2} {\Sb}_{\thetab} (\Xb_{ij})$ for all
$j=1, \ldots, n_i$ and $i=1,  \ldots, m$. \Black
Let  $\mathcal{G}_{\underline{h}}$ (${\underline{h}}:=(h_1, \ldots, h_m)$) be
the group of transformations of the form 
\begin{align*}\label{monotrans} 
 g_{h_i}: \Xb_{ij}  \mapsto  {g}_{h_i}(\Xb_{ij})
= h_i(\Xb_{ij}\pr \thetab) \thetab + \sqrt{1-
  (h_i(\Xb_{ij}\pr\thetab))^2} {\Sb}_{\thetab} (\Xb_{ij}), \quad
i=1,\ldots, m,
 \end{align*}
where the $h_i:[-1,1]\rightarrow [-1,1]$ are monotone continuous
nondecreasing functions such that $h_i(1)=1$ and $h_i(-1)=-1$. For any $m$-tuple of (possibly different) transformations $(g_{h_1}, \ldots,
g_{h_m}) \in \mathcal{G}_{\underline{h}}$, it is easy to verify that
$\| {g}_{h_i}(\Xb_{ij}) \|=1$; thus, $g_{h_i}$ is a monotone
transformation from ${\mathcal S}^{k-1}$ to ${\mathcal S}^{k-1}$,
$i=1, \ldots, m$. Note furthermore that  $g_{h_i}$ does not modify the signs
${\Sb}_{\thetab} (\Xb_{ij})$. Hence the
group of transformations $\mathcal{G}_{\underline{h}}$  is a
generating group for $\bigcup_{\underline{f}\in {\mathcal F}^m} {\rm
  P}\n_{\varthetab_0; {\underline{f}}}$  and 
the null is invariant under the action of
$\mathcal{G}_{\underline{h}}$. Letting  $R_{ij}$ denote the rank of
$\Xb_{ij}\pr\thetab$ among $\Xb_{i1}\pr\thetab, \ldots,
\Xb_{in_i}\pr\thetab$, $i=1, \ldots, m$,  it is now easy 
to conclude that the 
maximal invariant associated with $\mathcal{G}_{\underline{h}}$ is the
vector of signs ${\bf S}_{\thetab}(\Xb_{11}), \ldots,{\bf
  S}_{\thetab}(\Xb_{1n_1}),\ldots, {\bf S}_{\thetab}(\Xb_{m1}),\ldots,  {\bf
  S}_{\thetab}(\Xb_{mn_m})$ and ranks $R_{11}, \ldots, R_{1n_1},
\ldots, R_{m1}, \ldots, R_{mn_m}$. 

 As a consequence, we choose to base
our tests in this section on a rank-based version of the central
sequence $\Deltab\n_{\varthetab_0; \underline{f}}$, namely on
$$\utDelta_{\varthetab_0; \underline{K}}\n:=(
(\utDelta\n_{{\thetab}; {K}_{1}})', \ldots, (\utDelta\n_{{\thetab};
  {K}_{m}})^{\prime} )\pr$$ 
with 
$$\utDelta\n_{{\thetab}; {K}_{i}}=  n_i^{-1/2} \sum_{j=1}^{n_i} K_{i} \left(\frac{R_{ij}}{n_i+1} \right) \Sb_\thetab(\Xb_{ij}), \quad i=1, \ldots, m,$$
where $\underline{K}_{}:=(K_ {1}, \ldots,  K_{m})$ is a $m$-tuple of
\emph{score (generating) functions} satisfying

\

\noindent {\sc Assumption E}. The score functions $K_{i}$, $i=1,  \ldots, m$, are
continuous functions from $[0,1]$ to $\R$. 

\

\noindent The following result, which is a direct corollary (using
again the inner-sample independence and the mutual independence
between the $m$ samples) of Proposition~3.1 in Ley \emph{et
  al.}~(2013), characterizes the asymptotic behavior of
$\utDelta_{\varthetab_0; \underline{K}_{}}\n$ under any $m$-tuple of
densities with respective angular functions $g_1, \ldots, g_m$. 

\begin{Prop}  \label{asymplin}
Let Assumptions A, B, C and E  hold and consider $\underline{g}=(g_1, \ldots, g_m)\in\mathcal{F}^m$. Then the rank-based central sequence  $\utDelta_{\varthetab_0;\underline{K}_{}}\n$
\begin{itemize}
\item[(i)] is such that $ \utDelta_{\varthetab_0; \underline{K}_{}}\n- \Deltab_{\varthetab_0; \underline{K}_{};
    {\underline {g}}}\n=o_{\rm P}(1)$ under ${\rm P}\n_{\varthetab_0; \underline{g}}$ as $\ny$, where ($\tilde{G}_i$ standing
  for the common cdf of the $\Xb_{ij}\pr \thetab$'s under ${\rm
    P}\n_{\varthetab_0; {\underline{g}}}$, $i=1,\ldots, m$) 
$$\Deltab_{\varthetab_0; \underline{K}; {\underline{g}}}\n=((\Deltab\n_{{\thetab}; K_{1};
  g_1})', \ldots, (\Deltab\n_{{\thetab}; K_{m}; g_m})')\pr$$ 
with
$$\Deltab\n_{{\thetab}; K_{i}; g_i}:= n_i^{-1/2} \sum_{j=1}^{n_i}
K_{i}\left(\tilde{G}_i(\Xb_{ij}\pr \thetab)\right)
\Sb_\thetab(\Xb_{ij}), \quad i=1, \ldots, m.$$
In particular, for $\underline{K}=\underline{K}_{\underline{f}}:=(K_{f_1}, \ldots,  K_{f_m})$ with $K_{f_i}(u):=\varphi_{f_i}(\tilde{F}_i^{-1}(u)) (1-(\tilde{F}_i^{-1}(u))^2)^{1/2}$, $\utDelta_{\varthetab_0; \underline{K}_{\underline{f}}}\n$ is asymptotically equivalent to  the efficient  central sequence $\Deltab_{\varthetab_0; \underline{f}}\n$ under ${\rm P}\n_{\varthetab_0; {\underline {f}}}$.
\item[(ii)] is asymptotically normal under ${\rm P}\n_{\varthetab_0; \underline{g}}$ with mean zero and covariance matrix  $$\Gamb_{\varthetab_0; \underline{K}}:={\rm diag} \left( \frac{ {\mathcal J}_k({K}_{1})}{k-1}({\bf I}_k- \thetab \thetab\pr) , \ldots,  \frac{ {\mathcal J}_k({K}_{m})}{k-1}({\bf I}_k- \thetab \thetab\pr)\right),$$ \Black where ${\mathcal J}_k({K_{i}}):= \int_0^1 K_{i}^2(u) du$.
\item[(iii)] is asymptotically normal under ${\rm
    P}\n_{\varthetab_0+ n^{-1/2} \nub\n \tb\n;
    {\underline {g}}}$ ($\tb^{(n)}$ as
  in~(\ref{localt})) with mean $\Gamb_{\varthetab_0; \underline{K},
    {\underline g}} \tb$ ($\tb:= \lim_{\ny} \tb^{(n)}$) and
  covariance matrix 
\begin{equation*}\Gamb_{\varthetab_0;  \underline{K}, {\underline g}}:= {\rm diag} \left( \frac{ {\mathcal J}_k({K}_{1}, g_1)}{k-1}({\bf I}_k- \thetab \thetab\pr) , \ldots,\frac{ {\mathcal J}_k({K}_{m}, g_m)}{k-1}({\bf I}_k- \thetab \thetab\pr)\right),
\end{equation*}
where $\mathcal{J}_k(K_{i},g_i):= \int_{0}^{1} K_{i}(u) K_{g_i}(u) du$ for $i=1, \ldots, m$.
\item[(iv)] satisfies, under ${\rm P}\n_{\varthetab_0; {\underline {g}}}$ as $\ny$, the asymptotic linearity property
$$\utDelta_{\varthetab_0+n^{-1/2} \nub\n \tb\n;  \underline{K}}\n- \utDelta_{\varthetab_0;  \underline{K}}\n =- \Gamb_{\varthetab_0;  \underline{K}, {\underline {g}}}\tb\n + o_{\rm P}(1),$$
for $\tb\n=(\tb_1^{(n) \prime}, \ldots, \tb_m^{(n) \prime})\pr$ as in~(\ref{localt}). 
\end{itemize}
\end{Prop}

Similarly as for the pseudo-FvML test, our rank-based procedures are
not complete since we still need to estimate the common value
$\thetab$ of $\thetab_1, \ldots, \thetab_m$ under $\mathcal{H}_0$.  To
this end we will assume the existence of an estimator
$\hat{\varthetab}$  satisfying the following strengthened version of
Assumption D :

\

\noindent {\sc Assumption D'}. Besides $n^{1/2}\big(\nub\n\big)^{-1}$-consistency under ${\rm P}\n_{\varthetab_0;\underline{g}}$ for any $\underline{g}\in\mathcal{F}^m$, the estimator $\hat\varthetab\in (\mathcal{S}^{k-1})^m$ is further  \emph{locally and asymptotically discrete}, meaning that it only takes a bounded number of distinct values in $\varthetab_0$-centered balls of the form $\{
\mathbf{t}\in\R^{mk}   :  n^{1/2} \| \big(\nub\n\big)^{-1} \Black  (\mathbf{t}-{\varthetab _0}) \| \leq c\}$.\vspace{2mm}

\

Estimators satisfying the above assumption are  easy to
construct. Indeed the consistency is not a problem and  
the discretization condition  is a purely technical requirement (needed
to deal with these rank-based test statistics, see   pages~125 and 188
of Le~Cam and Yang~2000 for a discussion) with little practical
implications (in fixed-$n$ practice, such discretizations are
irrelevant as the radius can be taken arbitrarily large). We will
therefore tacitly assume that $\hat{\thetab} \in \mathcal{S}^{k-1}$
(and therefore $\hat\varthetab={\bf 1}_m \otimes \hat{\thetab}$) is
locally and asymptotically discrete throughout this section. Following
Lemma 4.4 in Kreiss (1987), the local discreteness allows to replace
in Part~(iv) of Proposition \ref{asymplin} non-random
perturbations of the form ${\varthetab} + n^{-1/2}\nub\n{\bf t}\n$ with ${\bf t}\n$ such that ${\varthetab} + n^{-1/2}\nub\n{\bf t}\n$ still belongs to ${\mathcal H}_0$ by a $n^{1/2} (\nub\n)^{-1}$-consistent estimator $\hat{\varthetab}:={\bf 1}_m \otimes \hat \thetab$.
  Based on the asymptotic result of Proposition \ref{asymplin} and letting 
  \begin{align*} \label{Mat2} & \Gamb_{\varthetab_0;   \underline{K},
      {\underline g}}^{\perp}  := \Gamb_{{\varthetab}_0;
      \underline{K}}^- \nonumber \\
&    \quad -  \Gamb_{{\varthetab}_0;
      \underline{K}}^-\Gamb_{{\varthetab}_0;  \underline{K},
      {\underline g}} \Umb_{\varthetab_0; {\nub}}\n [
    (\Umb_{\varthetab_0; {\nub}}\n)\pr \Gamb_{{\varthetab}_0;
      \underline{K}, {\underline g}}   \Gamb_{{\varthetab}_0;
      \underline{K}}^-  \Gamb_{{\varthetab}_0;  \underline{K},
      {\underline g}} \Umb_{\varthetab_0; {\nub}}\n]^-
    (\Umb_{\varthetab_0; {\nub}}\n)\pr\Gamb_{{\varthetab}_0;
      \underline{K}, {\underline g}}   \Gamb_{{\varthetab}_0;
      \underline{K}}^-,& & \end{align*}
%
the $\underline{g}$-valid rank-based test statistic 
we propose for the present ANOVA problem  corresponds to the quadratic form
$$Q_{ \underline{K}}({\underline g})\n:=(\utDelta\n_{\hat\varthetab;   \underline{K}})\pr {\Gamb}_{\hat{\varthetab};  \underline{K},  \underline{g}}^{\perp}\utDelta\n_{\hat \varthetab;  \underline{K}}.$$\vspace{0mm}

\noindent This test statistic  still depends on the
cross-information quantities  
\begin{equation}
  \label{eq:3}
  {\mathcal J}_k(K_{f_1},g_1),
\ldots, {\mathcal J}_k(K_{f_m},g_m)
\end{equation}
 and hence is only valid under fixed
$\underline{g}$. Therefore, exactly as for the pseudo-FvML tests of
the previous section, the
final step in our construction consists in estimating these quantities
consistently. For this define, for any $\rho\geq 0$, 
\begin{equation*} 
\tilde{\thetab}_i(\rho):= \hat{\thetab}+ n_{i}^{-1/2} \rho\, (k-1) ({\bf I}_k-\hat{\thetab}\hat{\thetab}\pr) \utDelta\n_{\hat{\thetab}; K_{i}}, \quad i=1, \ldots, m. 
\end{equation*}
Then, letting $\hat{\thetab}_i(\rho):=\tilde{\thetab}_i(\rho)/ \| \tilde{\thetab}_i(\rho) \|$, we consider the piecewise continuous quadratic form
$$
\rho \mapsto h_i\n(\rho):={\frac{k-1}{{\mathcal J}(K_{i})}}(\utDelta\n_{\hat{\thetab}; K_{i}})\pr \utDelta\n_{\hat{\thetab}_i(\rho); K_{i}}.
$$
Consistent estimators of the
quantities $  {\mathcal J}_k^{-1}(K_{1},g_1), \ldots, {\mathcal
  J}_k^{-1}(K_{m},g_m)$
 (and therefore readily of~\eqref{eq:3})  can
be obtained by taking  
\begin{equation*}
  \hat\rho_i:=\inf\{\rho>0:
h_i\n(\rho)<0 \}
\end{equation*}
for $i = 1, \ldots, m$ (see also   Ley \emph{et al.}~2013 for more details). Denoting by
$\hat{\mathcal J}_k(K_{i},g_i)$, for $i=1, \ldots, m$,
the resulting estimators, setting  $\hat{H}_{
  \underline{K},\underline{g}}:=\sum_{i=1}^m {r_i\n \hat{\mathcal
    J}_k^2({K}_{i}, g_i)}/{{\mathcal J}_k({K_{i}})}$ and letting ${\bf
  U}_{ij}:=K_{i} \left({\hat{R}_{ij}}/{(n_i+1)} \right) \Sb_{{\hat
    \thetab}}(\Xb_{ij})$, $i=1, \ldots, m $, ($\hat{R}_{ij}$ naturally
stands for the rank of $\Xb_{ij}\pr \hat{\thetab}$ among $\Xb_{i1}\pr
\hat{\thetab}, \ldots, \Xb_{in_i}\pr \hat{\thetab}$), the proposed
rank test $\utphi_{ \underline{K}}\n$ rejects the null hypothesis of
homogeneity of the locations when 
\begin{eqnarray*}
 \utQ_{ \underline{K}}\n &:=& (k-1)\sum_{i=1}^{m}  \frac{n_i}{{\mathcal J}_k(K_{i})} \bar{\Ub}_i\pr \bar{\Ub}_i  -(k-1) \hat{H}_{ \underline{K},\underline{g}}^{-1} \sum_{i,j=1}^{m} \frac{n_i n_{j}}{n} \frac{\hat{\mathcal J}_k(K_{i},g_i)}{{\mathcal J}_k(K_{i})} \frac{\hat{\mathcal J}_k(K_{j},g_{j})}{{\mathcal J}_k(K_{j})} \bar{\Ub}_i\pr \bar{\Ub}_{j} \\ 
\end{eqnarray*}
exceeds the $\alpha$-upper quantile of the chi-square distribution
with $(m-1)(k-1)$ degrees of freedom. This asymptotic behavior under
the null as well as the asymptotic distribution of $  \utQ_{
  \underline{K}}\n $ under a sequence of contiguous alternatives are
summarized in the following proposition. 

\begin{Prop} \label{ranktest} 
Let Assumptions~A, B, C and E hold and let $\hat\varthetab$ be an estimator such that Assumption D' holds. Then 
\begin{itemize}
\item[(i)] $ \utQ_{ \underline{K}}\n$ is asymptotically chi-square with $(m-1)(k-1)$ degrees of freedom  under  \linebreak $\bigcup_{\varthetab_0\in\mathcal{H}_0}\bigcup_{{\underline {g}} \in \mathcal{F}^{m}} \{ {\rm P}\n _{\varthetab_0;{\underline {g}}}\}$; 
\item[(ii)] $ \utQ_{ \underline{K}}\n$  is asymptotically non-central chi-square, still with $(m-1)(k-1)$ degrees of
freedom, but with non-centrality parameter~$$
 l_{\varthetab_0,{\bf t};  \underline{K}, {\underline {g}}}:= {\bf t}\pr\Gamb_{\varthetab_0;  \underline{K}, {\underline {g}}}  {\Gamb}_{\varthetab_0;  \underline{K},  \underline{g}}^{\perp}\Gamb_{\varthetab_0;  \underline{K}, {\underline {g}}} {\bf t}
$$
under ${\rm P}\n _{\varthetab_0+n^{-1/2} \nub\n \tb\n;{\underline {g}}}$, where $\tb^{(n)}$ is as in~(\ref{localt}) and $\tb:= \lim_{\ny} \tb^{(n)}$;
 \item[(iii)] the test  $\utphi_{ \underline{K}}\n$ which rejects the null hypothesis as soon as $ \utQ_{ \underline{K}}\n$ exceeds the $\alpha$-upper quantile of the chi-square distribution with $(m-1)(k-1)$ degrees of freedom has 
asymptotic level $\alpha$ under $\bigcup_{\varthetab_0\in {\mathcal H}_{0}}\bigcup_{ {\underline g}\in \mathcal{F}^{m}} \{ {\rm P}\n _{\varthetab_0; {\underline g}} \}$;
\item[(iv)] in particular, for $ \underline{K}=\underline{K}_{\underline{f}}:=(K_{f_1}, \ldots,  K_{f_m})$ with $K_{f_i}(u):=\varphi_{f_i}(\tilde{F}_i^{-1}(u)) (1-(\tilde{F}_i^{-1}(u))^2)^{1/2}$, $\utphi_{\underline{K}_{\underline f}}\n$  is locally and asymptotically most stringent,  
at asymptotic level~$\alpha$,  for $\bigcup_{\varthetab_0\in {\mathcal H}_{0}}
\bigcup_{{\underline g} \in \mathcal{F}^{m}} 
\{ {\rm P}\n _{\varthetab_0; {\underline  g}} \}$ against alternatives of the form
 $\bigcup_{\varthetab\notin{\mathcal H}_{0}} \{ {\rm
  P}\n _{\varthetab; {\underline  f}}\}$.
  \vspace{1mm}

\end{itemize}
\end{Prop}

\noindent Thanks to Proposition~\ref{asymplin}, the proof of this
result follows along the same lines as that of
Proposition~\ref{chideuxtwo} and is therefore omitted. 

We conclude this section by comparing the optimal pseudo-FvML test
$\phi\n$ with optimal rank-based tests $\utphi_{ \underline{K}_{\underline{f}}}\n$ for
several choices of $\underline{f}\in\mathcal{F}^m$ by means of
Pitman's asymptotic relative efficiency (ARE).  Letting
${\rm ARE}_{\varthetab_0;{\underline {g}}} (\phi_1\n,  \phi_2\n)$
denote the ARE of a test $\phi_1\n$ with respect to another test
$\phi_2\n$ under ${\rm P}\n_{\varthetab_0+n^{-1/2} \nub\n \tb\n; {\underline {g}}}$, we
have that
\begin{equation*}\label{ARE}
{\rm ARE}_{\varthetab_0;{\underline {g}}} (\utphi_{ \underline{K}_{\underline{f}}}\n,  \phi\n)= l_{\varthetab_0,{\bf t};  \underline{K}_{\underline{f}}, {\underline {g}}}/ l_{\varthetab_0,{\bf t}; \underline{\phi}, {\underline {g}}}.
\end{equation*}
In the homogeneous case ${\underline {g}}=(g_{1}, \ldots,  g_{1})$ (the angular density is the same for the $m$ samples) and if the same score function---namely, $K_{f_1}$---is used for the $m$ rankings (the test is therefore denoted by $\utphi_{{K_{f_1}}}\n$), the ratio in~(\ref{ARE}) simplifies into
\begin{equation}  {\rm
ARE}_{\varthetab_0;{\underline {g}}}(\utphi_{{K_{f_1}}}\n/ \phi\n)=\frac{{\cal J}_k^2(K_{f_1},g_{1})}{\mathcal{J}_k(K_{f_1}) D_{k,g_{1}}^2 B_{k,g_{1}}}. 
\label{arehomo} 
\end{equation}
Numerical values of the AREs in (\ref{arehomo}) are reported in Table \ref{AREtable1} for the three-dimensional setup under various angular densities and various choices of the score function $K_{f_1}$. More precisely, we consider the spherical linear, logarithmic and logistic distributions with respective angular functions
\begin{equation*}\label{linlog} f_{{\rm lin(a)}}(t):=t+a,\quad  f_{{\rm log(a)}}(t):=\log(t+a) \quad \mbox{  and } \end{equation*} \begin{equation*} \label{densities} f_{{\rm logis}(a,b)}(t):= \frac{a \; {\rm exp}(-b \; {\rm arccos}(t))}{(1+a \; {\rm exp}(-b \; {\rm arccos}(t)))^2}.\end{equation*}

\noindent The constants $a$ and $b$ are chosen so that all the above
functions are true angular functions satisfying Assumption~A. The
score functions associated with these angular functions are denoted by
$K_{\rm lin(a)}$ for $f_{{\rm lin(a)}}$, $K_{\rm log(a)}$ for $f_{{\rm
    log(a)}}$ and $K_{\rm logis (a,b)}$ for $f_{{\rm logis(a,b)}}$. For
the FvML distribution with concentration $\kappa$, the score function
will   be denoted by $K_{\phi_\kappa}$.

\begin{table}
\caption{Asymptotic relative efficiencies of (homogeneous) rank-based tests $\protect \utphi\n_{K_{f_1}}$ with respect to the pseudo-FvML test $\phi\n$ under various three-dimensional rotationally symmetric densities. \vspace{2mm}}
\centering
\fbox{
\begin{tabular}{|c|ccccccc|}
\hline
& \multicolumn{7}{|c|}{ARE($\utphi\n_{K_{f_1}}/ \phi\n$)} \\
\hline
\hline
 Underlying density & $\utphi\n_{K_{\phi_2}}$ &$\utphi\n_{K_{\phi_6}}$& $\utphi\n_{K_{\rm lin(2)}}$ & $\utphi\n_{K_{\rm lin(4)}}$  & $\utphi\n_{K_{\rm log(2.5)}}$& $\utphi\n_{K_{\rm logis(1,1)}}$ & $\utphi\n_{K_{\rm logis(2,1)}}$ \\
 \hline
 \hline
 {\rm FvML(1)} &  0.9744 &0.8787  & 0.9813 & 0.9979 & 0.9027  &0.9321 &  0.7364\\
 {\rm FvML(2)} & 1 & 0.9556  &0.9978 & 0.9586 & 0.9749  & 0.9823 &0.8480 \\
 {\rm FvML(6)} & 0.9555 & 1  &0.9381 &0.8517  &  0.9768 & 0.9911 &0.9280 \\
 {\rm Lin(2)} &  1.0539 &  0.9909  &1.0562 & 1.0215 & 1.0212 & 1.0247  &0.8796 \\
  {\rm Lin(4)} & 0.9709 & 0.8627  & 0.9795 &  1.0128  & 0.8856 &0.9231  &  0.7097 \\
 {\rm Log(2.5)} & 1.1610 & 1.1633 & 1.1514  &1.0413  & 1.1908 & 1.1625  &  1.0951 \\
  {\rm Log(4)} & 1.0182 & 0.9216 & 1.0261& 1.0347 & 0.9503 & 0.9741 &  0.7851 \\
   {\rm Logis(1,1)} & 1.0768  & 1.0865 & 1.0635 & 0.9991 & 1.0701  & 1.0962 & 0.9778 \\
    {\rm Logis(2,1)} &  1.3182 & 1.4426 & 1.2946 & 1.0893 & 1.4294 &1.3865 & 1.5544 \\
    \hline
\end{tabular}
}
\label{AREtable1}
\end{table}
 Inspection of Table \ref{AREtable1} confirms the theoretical results. As expected, the pseudo-FvML test $\phi\n$ dominates the rank-based tests under FvML densities, whereas rank-based tests mostly outperform the pseudo-FvML test under other densities, especially so when they are based on the score function associated with the underlying density (in which case the rank-based tests are optimal).

\section{Simulation results}\label{sec:Simul}

In this section, we perform a Monte Carlo study to compare the small-sample behavior of the pseudo-FvML test
$\phi^{(n)}$ and various rank-based tests $\utphi\n_{\underline{K}_{\underline f}}$ for the two-sample spherical location problem, that is, for an ANOVA with $m=2$. For this purpose, we generated $M=2,500$ replications of
four pairs  of mutually independent samples (with respective sizes $n_1=100$ and $n_2=150$) of $(k=)3$-dimensional rotationally symmetric random vectors 
$$
{\varepsb}_{\ell;ij_{i}}, \quad
\ell=1, 2, 3,4,  \ \ \ j_{i}=1,\ldots, n_{i}, \quad i=1,2,
$$
with FvML densities and linear densities: the ${\varepsb}_{1;1j_1}$'s have a FvML(15) distribution and the ${\varepsb}_{1;2j_2}$'s have a FvML(2) distribution; the ${\varepsb}_{2;1j_1}$'s have a Lin(2) distribution and the ${\varepsb}_{2;2j_2}$'s have a Lin(1.1) distribution; the ${\varepsb}_{3;1j_1}$'s have a FvML(15) distribution and the ${\varepsb}_{3;2j_2}$'s have a Lin(1.1) distribution and finally the ${\varepsb}_{4;1j_1}$'s have a Lin(2) distribution and the ${\varepsb}_{4;2j_2}$'s have a FvML(2) distribution.
 
 The rotationally symmetric vectors ${\varepsb}_{\ell;ij_{i}}$'s have all been generated with a common spherical location $\thetab_0=(\sqrt{3}/2, 1/2, 0)\pr$. Then,  each replication of the ${\varepsb}_{\ell;ij_i}$'s was transformed into
 \begin{equation*}\label{samples}
 \left\{ \begin{array}{l}
\Xb_{\ell;1j_{1}}={\varepsb}_{\ell;1j_{1}}, \quad \ell=1, 2, 3, 4, \ \ \ j_{1}=1,\ldots,n_{1} \\
\Xb_{\ell;2j_{2}; \xi}={\bf
O}_{\xi}{\varepsb}_{\ell;2j_{2}}, \quad \ell=1, 2, 3,4,  \ \ \ j_{2}=1,\ldots,n_{2},  \ \ \xi=0, 1, 2, 3,
\end{array} \right.
 \end{equation*}
where
$$
{\bf O}_{\xi}=\left(\begin{array}{ccc} \cos (\pi \xi / 16 \Black)  & - \sin (\pi \xi / 16 \Black) & 0
\\ \sin (\pi \xi / 16 \Black) & \cos (\pi \xi / 16 \Black) & 0 \\ 0&0&1 \end{array}\right).
$$

Clearly, the spherical locations of the $\Xb_{\ell;1j_{1}}$'s and the $\Xb_{\ell;2j_{2};
  0}$'s   coincide while the spherical location of the $\Xb_{\ell;2j_{2}; \xi}$'s,
$\xi=1,2,3$, is different from the spherical location of the $\Xb_{\ell;1j_{1}}$'s,
characterizing alternatives to the null hypothesis of common spherical
locations. Rejection frequencies based on the asymptotic chi-square
critical values at nominal level $5\%$ are reported in Table~\ref{simuresu} below. The inspection of the latter reveals expected results:

\begin{itemize}

\item[(i)] The pseudo-FvML test and all the rank-based tests are valid under heterogeneous densities. They reach the $5\%$ nominal level constraint under any considered pair of densities.

\item[(ii)] The comparison of the empirical powers reveals that when
  based on scores associated with the underlying distributions, the
  rank-based test performs   nicely. The pseudo-FvML test is clearly optimal in the FvML case. 
 
 \end{itemize}
 
\begin{table}
\caption{
Rejection frequencies (out of $M=2, 500$ replications), under the null
and under increasingly distant  alternatives, of the pseudo-FvML test $\phi^{(n)}$ and various rank-based tests $\protect \utphi_{(K_{\phi_{15}},K_{\phi_{2}})}\n$ (based on ${\rm FvML}(15)$ and ${\rm FvML}(2)$ scores), $\protect \utphi_{(K_{\rm Lin(2)}, K_{\rm Lin(1.1)})}\n$  (based on ${\rm Lin}(2)$ and ${\rm Lin}(1.1)$ scores), $\protect \utphi_{(K_{\rm Lin(2)}, K_{\phi_2})}\n$  (based on ${\rm Lin}(2)$ and ${\rm FvML}(2)$ scores) and $\protect \utphi_{(K_{\phi_{15}}, K_{\rm Lin(1.1)})}\n$  (based on ${\rm FvML}(15)$ and ${\rm Lin}(1.1)$ scores). Sample sizes are $n_1=100$ and $n_2=150$.  \vspace{2mm}
}
\centering
\fbox{
\begin{tabular}{|c|c|cccc|c|cccc|}
\hline
& & \multicolumn{4}{|c|}{{\small $ \xi $}} \\
\hline
Test & True densities & 0 & $1$ & $2$ & $3$  \\
\hline
$\phi\n$          &                             & .0592 &   .2684  &.8052 &   .9888    \\ 
$\utphi_{(K_{\phi_{15}},K_{\phi_{2}})}\n$   &        & .0696  &  .2952  &  .8276  &  .9900    \\ 
$\utphi_{(K_{\rm Lin(2)}, K_{\rm Lin(1.1)})}\n$    &   $(\phi_{15}, \phi_2)$     &  .0536 &   .2316  &  .7660   & .9756  \\
$\utphi_{(K_{\rm Lin(2)}, K_{\phi_2})}\n$    &                            &  .0656 &   .2952  &  .8160  &  .9894   \\
$\utphi_{(K_{\phi_{15}}, K_{\rm Lin(1.1)})}\n$      &        & .0544  &  .2308  &  .7716  &  .9772     \\[.5mm]
\hline
$\phi\n$      &    & .0480 &  .0596  &  .0792   & .1312    \\ 
$\utphi_{(K_{\phi_{15}},K_{\phi_{2}})}\n$  &   & .0472 &   .0568  &  .0948  &  .1340   \\ 
$\utphi_{(K_{\rm Lin(2)}, K_{\rm Lin(1.1)})}\n$    &   $({\rm Lin}({2}), {\rm Lin}({1.1}))$                         & .0464  &  .0604  &  .0892  &  .1424   \\
$\utphi_{(K_{\rm Lin(2)}, K_{\phi_2})}\n$             &               & .0520 &   .0588  &  .0920 &   .1440 \\
$\utphi_{(K_{\phi_{15}}, K_{\rm Lin(1.1)})}\n$  &   & .0480 &   .0580  &  .0856 &   .1340   \\[.5mm]
\hline
$\phi\n$     &    & .0508  &  .0684  &  .1044   & .1512    \\ 
$\utphi_{(K_{\phi_{15}},K_{\phi_{2}})}\n$   & & .0540 &   .0648  &  .1012  &  .1532   \\ 
$\utphi_{(K_{\rm Lin(2)}, K_{\rm Lin(1.1)})}\n$    &        $({\rm Lin}({2}), \phi_2)$                         & .0512 &  .0664  &  .1084  &  .1608   \\
$\utphi_{(K_{\rm Lin(2)}, K_{\phi_2})}\n$         &          & .0508 &   .0656  &  .1072 &   .1620 \\
$\utphi_{(K_{\phi_{15}}, K_{\rm Lin(1.1)})}\n$    &   &  .0496 &   .0628  &  .1004 &   .1516   \\[.5mm]
\hline
$\phi\n$          &                             & .0468 &   .1008  &.2908 &   .5760    \\ 
$\utphi_{(K_{\phi_{15}},K_{\phi_{2}})}\n$   &        & .0628  &  .1288  &  .3612 &  .6788   \\ 
$\utphi_{(K_{\rm Lin(2)}, K_{\rm Lin(1.1)})}\n$    &   $(\phi_{15},{\rm Lin}({1.1}))$     &  .0512 &   .1156  &  .3636   & .6892    \\
$\utphi_{(K_{\rm Lin(2)}, K_{\phi_2})}\n$    &                            &  .0616 &   .1220  &  .3620  &  .6768 \\
$\utphi_{(K_{\phi_{15}}, K_{\rm Lin(1.1)})}\n$      &        & .0504  &  .1180  &  .3660  &  .6916   \\[.5mm]
\hline

\end{tabular}
}
\label{simuresu}
\end{table}

\section{Real-data example}\label{realdat}

In this section, we evaluate the usefulness of our tests on a real-data example. The data consist of measurements of remanent
magnetization in red slits and claystones made at 2 different
locations in Eastern New South Wales, Australia. These data have
already been used in Embleton and McDonnell (1980). The rotationally
symmetric assumption in the two samples seems to be appropriate since
data are clearly concentrated. However, the specification of the
angular functions is not reasonable. 

The main question for the
practitioner is to test whether the remanent magnetization obtained in
those samples comes from a single source of magnetism or
not. Therefore, we test here the null hypothesis ${\mathcal H}_0:
\thetab_1=\thetab_2$ against
$\mathcal{H}_1:\thetab_1\neq\thetab_2$. For this purpose, we   used
the pseudo-FVML test  $\phi\n_{}$ and rank-based tests  
$\utphi_{({\rm lin}(1.1), {\rm FvML}(10) )}\n$ and $\utphi_{({\rm
    lin}(1.1), {\rm FvML}(100) )}\n$ based respectively on the couples
of linear and FvML scores $({\rm lin}(1.1), {\rm FvML}(10))$ and
$({\rm lin}(1.1), {\rm FvML}(100))$. The corresponding test statistics
are given by  
$$ Q_{{}}\n=5.96652,  \utQ_{({\rm lin}(1.1), {\rm
    FvML}(10))}\n=5.477525 \mbox{ and }   \utQ_{({\rm lin}(1.1), {\rm
    FvML}(100))}\n=5.525854.$$
 At the asymptotic nominal level $5\%$, the tests  $\phi\n_{}$,
$\utphi_{({\rm lin}(1.1), {\rm FvML}(10) )}\n$ and $\utphi_{({\rm lin}(1.1), {\rm FvML}(100) )}\n$ do not reject the null hypothesis of equality of the modal directions since the $5 \%$-upper quantile of the chi-square distribution with 2 degrees of freedom is equal to $5.991465$. 

\begin{figure}[htbp!] 
\centering
\includegraphics[height=7cm, width=10cm]{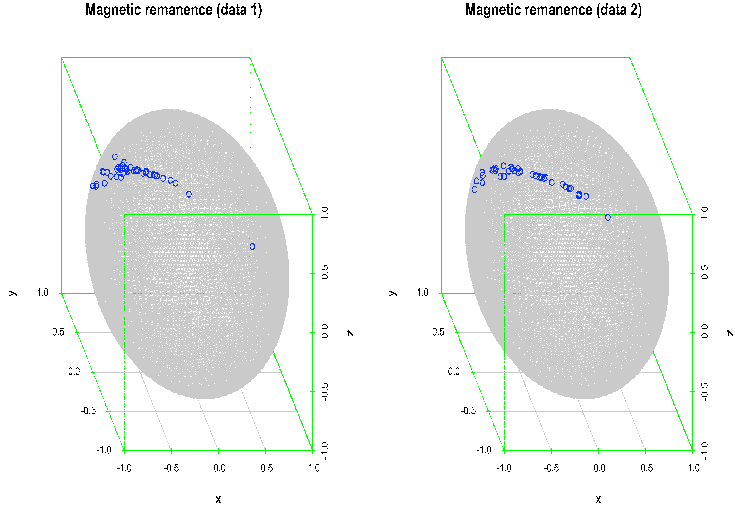}
\vspace{-0mm}
\caption{\small Measurements of remanent magnetization in red slits and claystones made at 2 different locations in Australia}
\label{data}
\end{figure} 

\noindent \Large\textbf{Appendix}\normalsize\vspace{0.6cm}

\setcounter{equation}{0}

\noindent {\bf Proof of Proposition \ref{lineasympbis}}. From Watson~(1983) (and the beginning of Section~\ref{rotsym}) we know that, under ${\rm P}\n_{\varthetab_0; \underline{g}}$, the sign vectors ${\bf S}_{\thetab}(\Xb_{ij})$ are independent of the scalar products $\Xb_{ij}\pr \thetab$, ${\rm E}_{g_i} \left[ {\bf S}_{\thetab}(\Xb_{ij})\right]=0$ and that
$${\rm E}_{g_i} \left[ {\bf S}_{\thetab}(\Xb_{ij}){\bf S}_{\thetab}\pr(\Xb_{ij}) \right]= \frac{1}{k-1}({\bf I}_k- \thetab \thetab\pr)$$
for $i=1, \ldots, m$ and for all $j=1,\ldots,n_i$. These results readily allow to obtain Part~(i) by applying the multivariate central limit theorem, while Part~(ii) 
 follows from the ULAN structure of the model in Proposition~\ref{ULAN} and Le Cam's third Lemma. \cqfd \vspace{4mm}

\noindent {\bf Proof of Proposition \ref{cross}}. We start by proving Part~(i). First note that easy computations yield (for $i=1, \ldots, m$)
\begin{eqnarray*}
\Deltab\n_{\hat\thetab; \phi_{\kappa_i}}&=& \kappa_in_i^{-1/2}\sum_{j=1}^{n_i}\left[\Xb_{i j}-(\Xb_{i j}\pr \hat\thetab) \hat\thetab\right]\\
 &=&\Deltab\n_{\thetab; \phi_{\kappa_i}}- \kappa_in_i^{-1/2}\sum_{j=1}^{n_i}\left[(\Xb_{i j}\pr \hat\thetab) \hat\thetab-(\Xb_{i j}\pr \thetab) \thetab\right] \\
 &=&\Deltab\n_{\thetab; \phi_{\kappa_i}}- {\bf V}_i\n- {\bf W}_i\n,
\end{eqnarray*}
where ${\bf V}_i\n:= \kappa_in_i^{-1} \sum_{j=1}^{n_i}\left[\Xb_{i j}\pr\thetab\right] n_i^{1/2} (\hat{\thetab}-\thetab)$ and ${\bf W}_i\n:=\hat{\thetab}\, \kappa_i n_i^{-1}(\sum_{j=1}^{n_i} \Xb_{i j}\pr) n_i^{1/2} (\hat{\thetab}-\thetab)$. Now, combining the delta method (recall that ${\bf I}_k -\thetab \thetab\pr$ is the Jacobian matrix of the mapping $h:\R^k\rightarrow\mathcal{S}^{k-1}:\xb \mapsto \frac{\xb}{\| \xb\|}$ evaluated at $\thetab$), the Law of Large Numbers and Slutsky's Lemma, we obtain that
\begin{eqnarray*}
 {\bf V}_i\n&=&\; \left(\kappa_in_i^{-1} \sum_{j=1}^{n_i}\Xb_{i j}\pr\thetab\right) n_i^{1/2} (\hat{\thetab}-\thetab) \\
&=& \;  \kappa_i{\rm E}_{g_i} [\Xb_{i j}\pr \thetab]  \; ({\bf I}_k- \thetab\thetab\pr) n_i^{1/2} (\hat{\thetab}-\thetab) + o_{\rm P}(1)\\
&=& \Gamb_{\thetab; g_i}^{\phi_{\kappa_i}} n_i^{1/2} (\hat{\thetab}-\thetab) + o_{\rm P}(1)
\end{eqnarray*}
under ${\rm P}\n_{\varthetab_0; \underline{g}}$ as $\ny$. Thus, the announced result follows as soon as we have shown that ${\bf W}_i\n$ is $o_{\rm P}(1)$ under ${\rm P}\n_{\varthetab_0; \underline{g}}$ as $\ny$. Using the same arguments as for ${\bf V}_i\n$, we have under ${\rm P}\n_{\varthetab_0; \underline{g}}$ and for $\ny$ that 
\begin{eqnarray*}
{\bf W}_i\n&=&\hat\thetab\, \left(\kappa_i n_i^{-1}\sum_{j=1}^{n_i} \Xb_{i j}\pr\right) n_i^{1/2} (\hat{\thetab}-\thetab) \\
&=& \hat\thetab\, \left(\kappa_in_i^{-1}\sum_{j=1}^{n_i} (\Xb_{i j}\pr)   \; ({\bf I}_k- \thetab\thetab\pr)\right)  n_i^{1/2} (\hat{\thetab}-\thetab) + o_{\rm P}(1)\\ 
  &=&   \hat\thetab \; \kappa_i{\rm E}_{g_i}\left[ \sqrt{1-(\Xb_{i j}\pr\thetab)^2}  ({\Sb}_{\thetab}(\Xb_{i j}))\pr \right] n_i^{1/2} (\hat{\thetab}-\thetab) + o_{\rm P}(1),
\end{eqnarray*}
which is $o_{\rm P}(1)$ from the boundedness of $\hat\thetab$ and since from Watson (1983) (see the proof of Proposition~\ref{lineasympbis} for more details) we know that 
$${\rm E}_{g_i}\left[ \sqrt{1-(\Xb_{i j}\pr\thetab)^2}  ({\Sb}_{\thetab}(\Xb_{i j}))\pr \right]={\rm E}_{g_i}\left[ \sqrt{1-(\Xb_{i j}\pr\thetab)^2} \right]{\rm E}_{g_i}\left[ ({\Sb}_{\thetab}(\Xb_{i j}))\pr \right]=\zerob\pr.$$
This concludes Part~(i) of the proposition. Regarding Part~(ii), let $\Xb$ be a random vector distributed according to an FvML distribution with concentration $\kappa$. Then, writing $c$ for the normalization constant, a simple integration by parts yields 
\begin{eqnarray*}
\label{egal}
C_{k, \phi_{\kappa}}= \kappa\,{\rm E}_{\phi_\kappa}[1-(\Xb\pr \thetab)^2] &=&  \kappa\, c \int_{-1}^1 (1-u^2) \; {\rm exp}(\kappa u)  (1-u^2)^{(k-3)/2}\; du \\
&=&  \kappa\, c \int_{-1}^1 \; {\rm exp}(\kappa u)  (1-u^2)^{(k-1)/2}\; du \\
 &=&c (k-1) \int_{-1}^1  u \; {\rm exp}(\kappa u)  (1-u^2)^{(k-3)/2}\; du \\
 &=& (k-1) \; {\rm E}_{\phi_\kappa}[\Xb\pr \thetab].
\end{eqnarray*}
The claim thus holds.  \cqfd

\vspace{4mm}

\noindent {\bf Proof of Proposition \ref{chideuxtwo}}. We start the
proof by showing that the replacement of $\thetab$ with
$\hat{\thetab}$ as well as the distinct estimators have no asymptotic
cost on $Q\n$. The consistency of $\hat{D}_{k, g_i}$, $\hat{E}_{k,
  g_i}$, $i=1,\ldots, m$, and
$\hat{H}_{\underline{\phi},\underline{g}}$ together with the $n^{1/2}
(\nub\n)^{-1}$-consistency of $\hat{\varthetab}$ entail that, using
Part~(i) of Proposition~\ref{cross}, 
$$Q\n=\left(\Deltab\n_{\varthetab_0;
    \underline{\phi}}-\Gamb_{\varthetab_0;
    \underline{g}}^{\underline{\phi}} \Umb\n \sqrt{n}
  \left(\hat\thetab -\thetab\right)\right)\pr {\Gamb}_{\varthetab_0;
  \underline{\phi},  \underline{g}}^{\perp}
\left(\Deltab\n_{\varthetab_0; \underline{\phi}}-\Gamb_{\varthetab_0;
    \underline{g}}^{\underline{\phi}} \Umb\n \sqrt{n}
  \left(\hat\thetab -\thetab\right)\right)+ o_{\rm P}(1)$$ 
under ${\rm P}\n_{\varthetab_0; \underline{g}}$ as $\ny$. Now,
standard algebra yields that $${\Gamb}_{{\thetab}; \underline{\phi},
  \underline{g}}^{\perp}\Gamb_{\varthetab_0;
  \underline{g}}^{\underline{\phi}} \Umb\n=(\Gamb_{\varthetab_0;
  \underline{g}}^{\underline{\phi}}\Umb\n)\pr{\Gamb}_{\varthetab_0;
  \underline{\phi}, \underline{g}}^{\perp}={\bf 0},$$ so that 
\begin{eqnarray*} Q\n &=& \left(\Deltab\n_{\varthetab_0;
      \underline{\phi}} \right)\pr {\Gamb}_{\varthetab_0;
    \underline{\phi},  \underline{g}}^{\perp} \Deltab\n_{\varthetab_0;
    \underline{\phi}}+ o_{\rm P}(1) \\ &=:& Q\n(\varthetab_0)+o_{\rm
    P}(1). 
\end{eqnarray*}
Both results from Proposition~\ref{lineasympbis} entail that since $\Gamb_{\varthetab_0; \underline{g}}^*
 {\Gamb}_{\varthetab_0; \underline{\phi}, \underline{g}}^{\perp}$ is  idempotent with trace $(m-1)(k-1)$,  $ Q\n(\thetab)$ (and therefore $ Q\n$) is asymptotically chi-square with $(m-1)(k-1)$ degrees of freedom
under ${\rm P}\n _{\varthetab_0; \underline{g}}$, and
asymptotically non-central chi-square, still with $(m-1)(k-1)$ degrees of freedom, and
with non-centrality parameter 
$\tb\pr \Gamb_{\varthetab_0; \underline{\phi}, \underline{g}} {\Gamb}_{{\thetab}; \underline{\phi},  \underline{g}}^{\perp} \Gamb_{\varthetab_0; \underline{\phi}, \underline{g}} \tb$ under ${\rm P}\n_{\varthetab_0+n^{-1/2} \nub\n \tb\n; \underline{g}}$. Parts (i) and (ii) follow.
Now, Part (iii) is a direct consequence of Part (i). Part (ii) of Proposition \ref{cross} and simple computations yield that $Q\n$ is asymptotically equivalent to the most stringent FvML test $Q_{\underline{\phi}}\n$ in (\ref{paramtest}). Part (iv) thus follows.  \vspace{4mm} \cqfd

\

\noindent {\large\bf Acknowledgements}\\

\noindent The research of Christophe Ley is supported by a Mandat de Charg\'e de Recherche from the Fonds National de la Recherche Scientifique, Communaut\' e fran\c caise de Belgique. \\

\

\bibliographystyle{elsarticle-harv}

\end{document}